\documentclass[a4paper,12pt,twoside]{article}
\usepackage[francais]{babel}
\usepackage{amssymb, amsmath, eucal}
\usepackage[all]{xy}
\usepackage{mathrsfs}
\begin{document}
\textwidth15.5cm
\textheight22.5cm
\voffset=-13mm
\newtheorem{The}{Theorem}[section]
\newtheorem{Lem}[The]{Lemma}
\newtheorem{Prop}[The]{Proposition}
\newtheorem{Cor}[The]{Corollary}
\newtheorem{Rem}[The]{Remark}
\newtheorem{Obs}[The]{Observation}
\newtheorem{Titre}[The]{\!\!\!\! }
\newtheorem{Conj}[The]{Conjecture}
\newtheorem{Question}[The]{Question}
\newtheorem{Prob}[The]{Problem}
\newtheorem{Def}[The]{Definition}
\newtheorem{Not}[The]{Notation}
\newtheorem{Claim}[The]{Claim}
\newtheorem{Ex}[The]{Example}
\newcommand{\C}{\mathbb{C}}
\newcommand{\R}{\mathbb{R}}
\newcommand{\N}{\mathbb{N}}
\newcommand{\Z}{\mathbb{Z}}
\newcommand{\Q}{\mathbb{Q}}

\begin{center}

{\Large\bf Stability of Strongly Gauduchon Manifolds under Modifications}

\end{center}

\begin{center}

 {\large Dan Popovici}

\end{center}

\vspace{1ex}

\noindent {\small {\bf Abstract.} In our previous works on deformation limits of projective and Moishezon manifolds, we introduced and made crucial use of the notion of {\it strongly Gauduchon} metrics as a reinforcement of the earlier notion of Gauduchon metrics. Using direct and inverse images of closed positive currents of type $(1, \, 1)$ and regularisation, we now show that compact complex manifolds carrying {\it strongly Gauduchon} metrics are stable under modifications. This stability property, known to fail for compact K\"ahler manifolds, mirrors the modification stability of balanced manifolds proved by Alessandrini and Bassanelli.}

\vspace{3ex}

\section{Introduction}\label{section:introd}

 Let $X$ be a compact complex manifold, $\mbox{dim}_{\C}X=n$. A Hermitian metric on $X$ will be identified throughout with the corresponding positive-definite $C^{\infty}$ $(1, \, 1)$-form $\omega$ on $X$. Recall that a Hermitian metric $\omega$ is said to be a {\it Gauduchon metric} if $\partial\bar{\partial}\omega^{n-1}=0$ on $X$ (cf. [Gau77]), a condition that can be reformulated as $\partial\omega^{n-1}$ being $\bar\partial$-closed. Gauduchon proved in [Gau77] that not only do these metrics always exist on any compact complex manifold $X$, but there is a Gauduchon metric in every conformal class of Hermitian metrics. 

 In [Pop09] we introduced the notion of a {\it strongly Gauduchon metric} (also referred to as an {\it sG metric} in the sequel) by requiring that, for a given Hermitian metric $\omega$, the above-mentioned $(n, \, n-1)$-form $\partial\omega^{n-1}$ be $\bar\partial$-exact on $X$ (rather than merely $\bar\partial$-closed). We showed that such a metric need not exist on an arbitrary $X$ and termed $X$ a {\it strongly Gauduchon manifold} if it carries a {\it strongly Gauduchon metric}. We went on to notice that when the $\partial\bar\partial$-lemma holds on $X$, the notions of Gauduchon and {\it strongly Gauduchon metrics} coincide, hence every such $X$ is a {\it strongly Gauduchon manifold}. This is because $d(\partial\omega^{n-1})=0$ if $\omega$ is a Gauduchon metric; since the pure-type form $\partial\omega^{n-1}$ is also $\partial$-exact in an obvious way, it must be $\bar\partial$-exact if the $\partial\bar\partial$-lemma holds on $X$. We then went on to characterise {\it strongly Gauduchon manifolds} starting from the following simple observation.

\begin{Lem}\label{Lem:ref-form} (Lemma 3.2. in [Pop09]) There exists a {\it strongly Gauduchon metric} on a given compact complex manifold $X$ of dimension $n$ if and only if there exists a $C^{\infty}$ $(2n-2)$-form $\Omega$ on $X$ such that: \\

\noindent $(a)$\, $\Omega=\overline{\Omega}$ (i.e. $\Omega$ is real); 

\noindent $(b)$\, $d\Omega=0$; 

\noindent $(c)$\, $\Omega^{n-1, \, n-1}>0$ on $X$ (i.e. the component of type $(n-1, \, n-1)$ of $\Omega$ is positive-definite).

\end{Lem}

 Indeed, if $\gamma$ is a {\it strongly Gauduchon metric} on $X$, set $\Omega^{n-1, \, n-1}:=\gamma^{n-1}$, take $\Omega^{n, \, n-2}$ to be any smooth $(n, \, n-2)$-form such that $\partial\gamma^{n-1}=-\bar\partial\Omega^{n, \, n-2}$ ($\Omega^{n, \, n-2}$ exists by the sG assumption on $\gamma$), set $\Omega^{n-2, \, n}:=\overline{\Omega^{n, \, n-2}}$ and $\Omega:=\Omega^{n, \, n-2} + \Omega^{n-1, \, n-1} + \Omega^{n-2, \, n}$. Conversely, given an $\Omega$ as in the above lemma, we use Michelsohn's procedure for extracting the $(n-1)^{st}$ root of a positive-definite $(n-1, \, n-1)$-form (cf. [Mic83]) and get a unique $C^{\infty}$ $(1, \, 1)$-form $\gamma$ on $X$ such that $\gamma^{n-1}= \Omega^{n-1, \, n-1}$. That $\gamma$ is an {\it sG metric} follows from the properties of $\Omega$. 

 As a consequence of this observation, we obtained the following intrinsic characterisation of {\it strongly Gauduchon manifolds}. \\

\begin{Prop}\label{Prop:intrinsic}(Proposition 3.3 in [Pop09]) A compact complex manifold $X$ is {\bf strongly Gauduchon} if and only if there is no non-zero current $T$ of type $(1, \, 1)$ on $X$ such that $T\geq 0$ and $T$ is $d$-exact on $X$. 

\end{Prop}

 The object of the present work is to show that the {\it strongly Gauduchon} property of compact complex manifolds is stable under modifications (i.e. proper, holomorphic, bimeromorphic maps). This provides a sharp contrast to the K\"ahlerness of these manifolds which is only preserved under blowing up (smooth) submanifolds ([Bla58]).

\begin{The}\label{The:mod-stab} Let $\mu:\widetilde{X}\rightarrow X$ be a modification of compact complex manifolds $X$ and $\widetilde{X}$. Then $\widetilde{X}$ is a {\bf strongly Gauduchon manifold} if and only if $X$ is a {\bf strongly Gauduchon manifold}.

\end{The}

 This result parallels the main result of Alessandrini and Bassanelli in [AB95] (see also [AB91] and [AB93]) which asserts that {\it balanced manifolds} enjoy the same stability property under modifications as above. Recall that {\it balanced manifolds} were given in [Mic83] an intrinsic characterisation in terms of non-existence of non-trivial positive currents of bidegree $(1, \, 1)$ that are components of a boundary. Our criterion listed as Proposition \ref{Prop:intrinsic} above is the analogous intrinsic characterisation of the weaker notion of {\it strongly Gauduchon manifolds}. The proof of Theorem \ref{The:mod-stab} will draw on those of the main results in [AB91], [AB93] and [AB95] with certain arguments handled slightly differently while others are considerably simplified by the fact that $d$-closed positive $(1, \, 1)$-currents always admit unambiguously defined inverse images constructed from their local potentials, unlike the much more delicate-to-handle $\partial\bar\partial$-closed positive $(1, \, 1)$-currents that were relevant to the case of {\it balanced manifolds}. Inverse images for this latter class of currents were painstakingly constructed in [AB93] and a unique choice was shown to enjoy the necessary cohomological properties, rendering the case treated in [AB93] and [AB95] conspicuously more involved than ours.

 The motivation for proving stability properties of {\it strongly Gauduchon manifolds} stems from the relevance of this notion to the study of limits of compact complex manifolds under holomorphic deformations. It has played a key role in proving that the deformation limit of projective (or merely Moishezon) manifolds is Moishezon (cf. [Pop09], resp. [Pop10]) and has similar striking consequences in efforts to tackle analogous but more general conjectures on deformations of {\it class} ${\cal C}$ manifolds. Besides their modification stability exhibited in the present text, {\it strongly Gauduchon manifolds} are likely to evince stability properties under holomorphic deformations generalising those shown in [Pop09] and [Pop10], a study of which is well underway and will hopefully be the subject of a future paper.

\section{\bf Proof of Theorem \ref{The:mod-stab}}

 Let $\mu :\widetilde{X}\rightarrow X$ be a modification of compact complex manifolds and denote $n=\mbox{dim}_{\C}\widetilde{X}=\mbox{dim}_{\C}X$. Let $E$ be the exceptional divisor of $\mu$ on $\widetilde{X}$ and let $S\subset X$ be the analytic subset such that the restriction $\mu_{|\widetilde{X}\setminus E}:\widetilde{X}\setminus E\longrightarrow X\setminus S$ is a biholomorphism. Theorem \ref{The:mod-stab} falls into two parts.

\begin{The}\label{The:downsG-upsG} If $\mu :\widetilde{X}\rightarrow X$ is a modification of compact complex manifolds and $X$ is {\bf strongly Gauduchon}, then $\widetilde{X}$ is again {\bf strongly Gauduchon}.

\end{The}

\noindent {\it Proof.} We proceed by contradiction. Suppose that $\widetilde{X}$ is not {\it strongly Gauduchon}. Then, by Proposition \ref{Prop:intrinsic}, there exists a current $T\neq 0$ of type $(1, \, 1)$ on $\widetilde{X}$ such that

$$T\geq 0 \hspace{3ex} \mbox{and} \hspace{3ex} T\in\mbox{Im}\, d \hspace{2ex} \mbox{on}\hspace{2ex} \widetilde{X}.$$

 By compactness of $\widetilde{X}$, the map $\mu$ is proper and therefore the direct image under $\mu$ of any current on $\widetilde{X}$ is well-defined. Thus $\mu_{\star}T$ is a well-defined current of type $(1, \, 1)$ on $X$. It is clear that

$$\mu_{\star}T\geq 0 \hspace{3ex} \mbox{and} \hspace{3ex} \mu_{\star}T\in\mbox{Im}\, d \hspace{2ex} \mbox{on}\hspace{2ex} X.$$

\noindent Indeed, for every $C^{\infty}$ $(1, \, 1)$-form $\omega >0$ on $X$, we have

$$\int\limits_X\mu_{\star}T\wedge\omega^{n-1}=\int\limits_XT\wedge(\mu^{\star}\omega)^{n-1}\geq 0,$$

\noindent a fact that proves the positivity of $\mu_{\star}T$, while the $d$-exactness follows from $\mu_{\star}$ commuting with $d$. Now we have the following dichotomy.

 If $\mu_{\star}T$ is non-zero, we get a contradiction to the {\it strongly Gauduchon} assumption on $X$ thanks to Proposition \ref{Prop:intrinsic}.  

 If $\mu_{\star}T=0$ on $X$, we show that $T=0$ on $\widetilde{X}$, contradicting the choice of $T$. Indeed, if $\mu_{\star}T=0$, the support of $T$ must be contained in the support of $E$. Since $T$ is a closed positive current of bidegree $(1, \, 1)$ and the irreducible components $E_j$ of $E$ are all of codimension $1$ in $\widetilde{X}$, a classical theorem of support (see e.g. [Dem97, Chapter III, Corollary 2.14]) forces $T$ to have the shape

$$T=\sum\limits_{j\in J}\lambda_j\, [E_j], \hspace{3ex} \mbox{with coefficients} \hspace{2ex} \lambda_j\geq 0 \hspace{2ex} \mbox{and some index set} \hspace{1ex} J.$$

 Now there are two cases. If all the irreducible components of $S$ are of codimension $\geq 2$ in $X$, then $\mbox{codim}_X\mu(E_j)\geq 2$ for every $j\in J$. All we have to do is repeat the argument of [AB91, p. 5] that we now recall for the reader's convenience. By [GR70, p. 286], for every $i\geq 0$, there exists a vector subspace $H_i^{\star}(E)\subset H_i(E)$ and a commutative diagram whose rows are short exact sequences featuring the homology groups $H_i$ of the various spaces involved:

$$\begin{array}{ccccccccc}\nonumber 0 & \longrightarrow & H_i^{\star}(E) & \hookrightarrow & H_i(E) & \stackrel{\beta_i}{\longrightarrow} & H_i(S) & \rightarrow & 0\\
\nonumber & &\Vert & &\downarrow & & \downarrow & & \\
\nonumber 0 & \longrightarrow & H_i^{\star}(E) & \longrightarrow & H_i(\widetilde{X}) & \stackrel{\alpha_i}{\longrightarrow} & H_i(X) & \rightarrow & 0,\end{array}$$

\noindent where $\hookrightarrow$ stands for inclusion. If we denote by $\{\,\,\}_E$ (respectively $\{\,\,\}_{\widetilde{X}}$) the homology class of a $d$-closed current of dimension $2(n-1)$ in the ambient space $\mbox{Supp\,E}$ (respectively ${\widetilde{X}}$), we see that 

$$\beta_{2(n-1)}\{T\}_E=\sum\limits_j\lambda_j\, \beta_{2(n-1)}\{[E_j]\}_E = 0$$

\noindent since the direct images under $\mu$ of the currents of integration $[E_j]$ are closed positive $(1, \, 1)$-currents on $X$ supported on the analytic subset $S$ with $\mbox{codim}_XS\geq 2$, hence they must vanish by another classical theorem of support (see e.g. [Dem97, Chapter III, Corollary 2.11]). Thus, from the top exact sequence, we get that $\{T\}_E$ belongs to $H_{2(n-1)}^{\star}(E)$. The diagram being commutative, the image of $\{T\}_E\in H_{2(n-1)}^{\star}(E)$ in $H_{2(n-1)}(\widetilde{X})$ under the injective arrow of the bottom exact sequence is $\{T\}_{\widetilde{X}}$. Meanwhile $\{T\}_{\widetilde{X}}=0$ since $T$ is $d$-exact on $\widetilde{X}$. We get that $\{T\}_E=0$. This means that the restriction $T_{|\mbox{Supp\, E}}$ is a $d$-exact current of bidegree $(0, \, 0)$ on $\mbox{Supp}\, E$ (since it is of the same bidimension $(n-1, \, n-1)$ as that of $T$ on $\widetilde{X}$). Since the only $d$-exact current of type $(0, \, 0)$ is the zero current, we must have $\lambda_j=0$ for every $j$. Hence $T=0$ as a current on $\widetilde{X}$, a contradiction.

 If $S$ has irreducible components $S_j$, $j\in J_0$, of codimension $1$ in $X$, then

 $$\mu^{-1}(S_j)=E_j \hspace{2ex} \mbox{and} \hspace{2ex} \mu^{\star}([S_j])=[E_j], \hspace{3ex}  j\in J_0\subset J.$$

\noindent Thus $T=\sum\limits_{j\in J_0}\lambda_j\, [E_j] + \sum\limits_{j\in J\setminus J_0}\lambda_j\, [E_j]$ and $0=\mu_{\star}T = \sum\limits_{j\in J_0}\lambda_j\, [S_j]$. Indeed, $\mu_{\star}[E_j]=[S_j]$ for all $j\in J_0$, while $\mu_{\star}[E_j]=0$ for every $j\in J\setminus J_0$ since it is a closed positive $(1, \, 1)$-current whose support is contained in an analytic subset of codimension $\geq 2$ in $X$. Hence $\lambda_j=0$ for all $j\in J_0$. Consequently, $T = \sum\limits_{j\in J\setminus J_0}\lambda_j\, [E_j]$ and we are now in the previous case where we showed that $T=0$, a contradiction. The proof is complete.  \hfill $\Box$

\vspace{2ex}

 We now prove the reverse statement.

\begin{The}\label{The:upsG-downsG} If $\mu :\widetilde{X}\rightarrow X$ is a modification of compact complex manifolds and $\widetilde{X}$ is {\bf strongly Gauduchon}, then $X$ is again {\bf strongly Gauduchon}.

\end{The}

\noindent {\it Proof.} We proceed once more by contradiction. Suppose that $X$ is not {\it strongly Gauduchon}. Then, in view of Proposition \ref{Prop:intrinsic}, there exists a current $T\neq 0$ of type $(1, \, 1)$ on $X$ such that

$$T\geq 0 \hspace{3ex} \mbox{and} \hspace{3ex} T= dS \hspace{2ex} \mbox{for some real} \hspace{1ex} 1\!\!-\!\mbox{current} \hspace{1ex} S \hspace{1ex} \mbox{on}\hspace{2ex} X.$$

 We shall show that the inverse image current $\mu^{\star}T$ is a well-defined $(1, \, 1)$-current on $\widetilde{X}$ enjoying the same properties as $T$ on $X$, thus contradicting the {\it strongly Gauduchon} assumption on $\widetilde{X}$ in view of Proposition \ref{Prop:intrinsic}.

 Although the inverse image of an arbitrary current is not defined in general, the inverse image of a $d$-closed positive $(1, \, 1)$-current is well-defined under $\mu$ by the inverse images of its local $\partial\bar\partial$-potentials (see e.g. [Meo96]). Indeed, following [Meo96], for every open subset $U\subset X$ such that $T_{|U}=i\partial\bar\partial\varphi$ for a psh function $\varphi$ on $U$, one defines $(\mu^{\star}T)_{|\mu^{-1}(U)}:=i\partial\bar\partial(\varphi\circ\mu)$. The psh function $\varphi\circ\mu$ is $\not\equiv -\infty$ on every connected component of $\mu^{-1}(U)$ since $\mu$ has generically maximal rank and the local pieces $(\mu^{\star}T)_{|\mu^{-1}(U)}$ glue together into a globally defined $d$-closed positive $(1, \, 1)$-current $\mu^{\star}T$ on $\widetilde{X}$ that is independent of the choice of open subsets $U\subset X$ and local potentials $\varphi$.

 It is clear that $\mu^{\star}T$ is not the zero current on $\widetilde{X}$. Indeed, if we had $\mu^{\star}T=0$, the support of $T$ would be contained in $S$. If all the irreducible components of $S$ were of codimension $\geq 2$ in $X$, a classical theorem of support (see e.g. [Dem97, Chapter III, Corollary 2.11]) would guarantee that the closed positive $(1, \, 1)$-current $T$ must be the zero current on $X$, a contradiction. If $S$ had certain global irreducible components $S_j$ of codimension $1$ in $X$, another theorem of support (cf. [Dem97, Chapter III, Corollary 2.14]) would ensure that $T$ has the shape $T=\sum\lambda_j\, [S_j]$ for some constants $\lambda_j\geq 0$. Then $\mu^{\star}[S_j]$ would be the current of integration on the inverse-image divisor $\mu^{-1}(S_j)\subset\widetilde{X}$ and $\mu^{\star}T$ cannot be the zero current unless $\lambda_j=0$ for all $j$. However, in this event $T=0$ on $X$, a contradiction.

 The only thing that has yet to be checked before reaching the desired contradiction is that the non-trivial $d$-closed positive $(1, \, 1)$-current $\mu^{\star}T$ is $d$-exact on $\widetilde{X}$. Since the $1$-current $S$ cannot be pulled back to $\widetilde{X}$ (the local potential technique is no longer available), we shall use Demailly's regularisation-of-currents theorem [Dem92, Main Theorem 1.1] to get a sequence $(v_j)_{j\in\N}$ of $C^{\infty}$ $(1, \, 1)$-forms on $X$ such that every $v_j$ lies in the same Bott-Chern (hence also De Rham) cohomology class as $T$ with convergence

$$v_j\longrightarrow T \hspace{2ex} \mbox{weakly as} \hspace{2ex} j\rightarrow +\infty, \hspace{2ex} \mbox{while}  \hspace{2ex} v_j\geq -C\omega, \hspace{2ex} j\in\N,$$

\noindent where $\omega$ is any Hermitian metric on $X$ fixed beforehand and $C>0$ is a constant independent of $j\in\N$.

 Since $T$ is $d$-exact and cohomologous to each $v_j$, every form $v_j$ is $d$-exact. Thus, for all $j\in\N$, $v_j=du_j$ for some $C^{\infty}$ $1$-form $u_j$ on $X$. Unlike $S$, the $C^{\infty}$ forms $u_j$ have inverse images under $\mu$ and we get

\begin{equation}\label{eqn:up-conv}\mu^{\star}v_j=d(\mu^{\star}u_j)\longrightarrow \mu^{\star}T \hspace{3ex} \mbox{weakly as} \hspace{2ex} j\rightarrow +\infty,\end{equation}

\noindent after possibly extracting a subsequence.  Indeed, it was shown in [Meo96, Proposition 1] that for every sequence of $d$-closed positive $(1, \, 1)$-currents $T_j$ converging weakly to $T$, the sequence of inverse-image currents $\mu^{\star}T_j$ converges weakly to $\mu^{\star}T$. In our case, the $(1, \, 1)$-forms $v_j$ are not necesarily {\it positive} but only {\it almost positive} (the negative part being uniformly bounded by $C\omega$). We now spell out the reason why $\mu^{\star}v_j$ converges weakly to the current $\mu^{\star}T$ in this slightly more general context. The argument is virtually the same as that of [Meo96].

 Pick any $C^{\infty}$ $(1, \, 1)$-form $\alpha$ in the Bott-Chern class of the forms $v_j$ ($=$ the class of $T$). Then, for every $j\in\N$, we have

$$v_j=\alpha + i\partial\bar\partial\psi_j\geq -C\omega \hspace{3ex} \mbox{on} \hspace{2ex} X, $$

\noindent with $C^{\infty}$ functions $\psi_j : X\rightarrow \R$ that we normalise such that $\int\limits_X\psi_j\,\omega^n=0$ for every $j$. This normalisation makes $\psi_j$ unique. Applying the trace w.r.t. $\omega$ and using the corresponding Lapacian $\Delta_{\omega}(\cdot)=\mbox{Trace}_{\omega}(i\partial\bar\partial(\cdot))$, we get

$$\Delta_{\omega}\psi_j=\mbox{Trace}_{\omega}(v_j-\alpha), \hspace{2ex} j\in\N.$$

 Applying now the Green operator $G$ of $\Delta_{\omega}$ and using the normalisation of $\psi_j$, we get

$$\psi_j=G\,\mbox{Trace}_{\omega}(v_j-\alpha), \hspace{2ex} j\in\N.  $$

 Since $G$ is a compact operator from the Banach space of bounded Borel measures on $X$ to $L^1(X)$ and since the forms $v_j$ converge weakly to $T$, we infer that some subsequence $(\psi_{j_k})_k$ converges to a limit $\psi\in L^1(X)$ in $L^1(X)$-topology. Thus the weak continuity of $\partial\bar\partial$ gives

$$T=\lim\limits_k(\alpha + i\partial\bar\partial\psi_{j_k}) = \alpha + i\partial\bar\partial\psi  \hspace{3ex} \mbox{on} \hspace{2ex} X.$$

 Now the sequence $(\psi_j)_j$ is uniformly bounded above on $X$ by some constant $C_1>0$ thanks to the normalisation imposed on $\psi_j$ and the Green-Riesz representation formula for $\psi_j$, $\Delta_{\omega}$ and $G$. Hence the sequence $(\psi_j\circ\mu)_j$ is uniformly bounded above on $\widetilde{X}$ by $C_1>0$. On the other hand, $\psi_{j_k}\circ\mu$ converges almost everywhere to $\psi\circ\mu$ on $\widetilde{X}$. Since the forms $i\partial\bar\partial(\psi_{j_k}\circ\mu)$ are uniformly bounded below on $\widetilde{X}$ by $-(\mu^{\star}\alpha + C\mu^{\star}\omega)$, the almost psh functions $\psi_{j_k}\circ\mu$ can be simultaneously made psh on small open subsets of $\widetilde{X}$ by the addition of a same locally defined smooth psh function. We can thus apply the classical result stating that a sequence of psh functions that are locally uniformly bounded above either converges locally uniformly to $-\infty$ (a case that is ruled out in our present situation), or has a subsequence that converges in $L^1_{loc}$ topology (see e.g. [Hor94, Theorem 3.2.12., p. 149]). We infer that the almost psh functions $\psi_{j_k}\circ\mu$ actually converge in $L^1(\widetilde{X})$-topology (hence also in the weak topology of distributions) and implicitly the forms

$$\mu^{\star}v_{j_k}=\mu^{\star}\alpha + i\partial\bar\partial(\psi_{j_k}\circ\mu)$$

\noindent converge weakly to the current $\mu^{\star}T=\mu^{\star}\alpha + i\partial\bar\partial(\psi\circ\mu)$. Thus the convergence statement (\ref{eqn:up-conv}) is proved.

 Since the De Rham class is continuous w.r.t. the weak topology of currents and since each form $\mu^{\star}v_j=d(\mu^{\star}u_j)$ has vanishing De Rham class, the limit current $\mu^{\star}T$ must have vanishing De Rham class. Equivalently, $\mu^{\star}T$ is $d$-exact, providing a contradiction to the {\it strongly Gauduchon} assumption on $\widetilde{X}$ in view of Proposition \ref{Prop:intrinsic}. The proof is complete. \hfill $\Box$

\vspace{6ex}

\noindent {\bf References.} \\

\noindent [AB91]\, L. Alessandrini, G. Bassanelli --- {\it Smooth Proper Modifications of Compact K\"ahler Manifolds} --- Proc. Internat. Workshop on Complex Analysis (Wuppertal 1990); Complex Analysis, Aspects of mathematics, {\bf E17}, Vieweg, Braunschweig (1991), 1-7.

\vspace{1ex}

\noindent [AB93]\, L. Alessandrini, G. Bassanelli --- {\it Metric Properties of Manifolds Bimeromorphic to Compact K\"ahler Spaces} --- J. Diff. Geom. {\bf 37} (1993), 95-121.

\vspace{1ex}

\noindent [AB95]\, L. Alessandrini, G. Bassanelli --- {\it Modifications of Compact Balanced Manifolds} --- C. R. Acad. Sci. Paris, t {\bf 320}, S\'erie I (1995), 1517-1522. 

\vspace{1ex}

\noindent [Bla58]\, A. Blanchard --- {Les vari\'et\'es analytiques complexes} --- Ann. Sci. \'Ecole Norm. Sup. {\bf 73} (1958), 157-202.

\vspace{1ex}

\noindent [Dem92] \, J. -P. Demailly --- {\it Regularization of Closed Positive Currents and Intersection Theory} --- J. Alg. Geom., {\bf 1} (1992), 361-409.

\vspace{1ex}

\noindent [Dem 97] \, J.-P. Demailly --- {\it Complex Analytic and Algebraic Geometry}---http://www-fourier.ujf-grenoble.fr/~demailly/books.html

\vspace{1ex}

\noindent [Gau77]\, P. Gauduchon --- {\it Le th\'eor\`eme de l'excentricit\'e nulle} --- C.R. Acad. Sc. Paris, S\'erie A, t. {\bf 285} (1977), 387-390.

\vspace{1ex}

\noindent [GR70]\, H. Grauert, O. Riemenschneider --- {\it Verschwindungss\"atze f\"ur analytische Kohomologiegruppen auf komplexen R\"aumen} --- Invent. Math. {\bf 11} (1970), 263-292.

\vspace{1ex}

\noindent [Hor94]\, L. H\"ormander --- {\it Notions of Convexity} --- Progress in Mathematics, v. 127, Birkh\"auser, Boston 1994.

\vspace{1ex}

\noindent [Meo96]\, M. Meo --- {\it Image inverse d'un courant positif ferm\'e 
par une application analytique surjective} --- C. R. Acad. Sci. Paris, t. {\bf 322}, S\'erie I, p. 1141 - 1144, 1996.

\vspace{1ex}

\noindent [Mic83]\, M. L. Michelsohn --- {\it On the Existence of Special Metrics in Complex Geometry} --- Acta Math. {\bf 143} (1983) 261-295.

\vspace{1ex}

\noindent [Pop09]\, D. Popovici --- {\it Limits of Projective Manifolds under Holomorphic Deformations} --- arXiv e-print math.AG/0910.2032v1.

\vspace{1ex}

\noindent [Pop10]\, D. Popovici --- {\it Limits of Moishezon Manifolds under Holomorphic Deformations} --- arXiv e-print math.AG/1003.3605v1.

\vspace{6ex}

\noindent Universit\'e Paul Sabatier

\noindent Institut de Math\'ematiques de Toulouse

\noindent 118, route de Narbonne

\noindent 31 062, Toulouse, France

\noindent Email: popovici@math.ups-tlse.fr

\end{document}